%% file: cross_matching.tex
\documentclass[10pt,twocolumn,letterpaper]{article}

\usepackage[pagenumbers]{cvpr} 

\input{preamble}

\definecolor{cvprblue}{rgb}{0.21,0.49,0.74}
\usepackage[pagebackref,breaklinks,colorlinks,allcolors=cvprblue]{hyperref}

\def\paperID{10483} 
\def\confName{CVPR}
\def\confYear{2025}

\title{Cross-Matched Interval Prevalence of High Dimensional Point Clouds}

\author{Jonathan M. Mousley\\
Duke University\\
Department of Mathematics\\
{\tt\small jonathan.mousley@duke.edu}
\and
Paul Bendich\\
Geometric Data Analytics\\
Duke University\\
Department of Mathematics\\
{\tt\small bendich@math.duke.edu}
}

\begin{document}
\maketitle

{
    \small
    \bibliographystyle{ieeenat_fullname}
    \bibliography{main}
}

\input{sec/0_abstract}
\input{sec/1_intro}
\input{sec/2_prelim}
\input{sec/3_cm}
\input{sec/4_implementation}
\input{sec/5_param}
\input{sec/6_exp}
\input{sec/7_conclusion}

\section*{Acknowledgements}
Research by the first author was partially funded by NSF DGE 2139754.
Research by the second author was partially funded by the National Institute of Aerospace (NIA) under sub-award C21-202066-GDA. We are grateful to Erin Taylor for helpful discussions about graph clustering.


\end{document}

%% file: preamble.tex
\usepackage{tikz,tikz-cd}
\usetikzlibrary{positioning, arrows.meta}
\newcommand{\dgm}{\text{dgm}}
\newcommand{\barcode}{\text{bar}}
\newcommand{\affinity}{\text{aff}}
\usepackage{amsfonts}
\usepackage{graphicx}
\usepackage{algorithm}
\usepackage{algpseudocode}
\usepackage{amsthm}
\usepackage{caption}
\theoremstyle{definition}
\newtheorem{definition}{Definition}[section]
\theoremstyle{theorem}

\theoremstyle{remark}
\newtheorem{remark}{Remark}[section]
\newcommand{\prev}{{\text{prev}}}
\usepackage{bm}
\graphicspath{ {./images/} }

%% file: sec/0_abstract.tex
\begin{abstract}
Topological Data Analysis (TDA) has been applied with success to solve problems across many scientific disciplines. However, in the setting of a point cloud $X$ sampled from a shape $\mathcal{S}$ of low intrinsic dimension embedded within high ambient dimension $\mathbb{R}^D$, persistent homology, a key element to many TDA pipelines, suffers from two problems. First, when relatively small amounts of noise are introduced to the point cloud, persistent homology is unable to recover the true shape of $\mathcal{S}$. Secondly, the computational complexity of persistent homology scales poorly with the size of a point cloud. Although there is recent work that addresses the first issue via topological bootstrapping methods and topological prevalence, these new techniques still fall victim to the second issue. Here we introduce the cross-matched prevalence image (CMPI), an image which approximates the topological prevalent information of said point cloud, requiring only computations of persistent homology on the scale of samples of the point cloud and not the entire point cloud itself. We compute the CMPI for high dimensional synthetic data, demonstrating that it performs similarly in noise robustness experiments and accurately captures prevalent topological features as compared to previous topological bootstrapping methods.
\end{abstract}

%% file: sec/1_intro.tex
\section{Introduction}
\label{sec:intro}
Topological data analysis (TDA) has evolved along several different trajectories \cite{ELZ,Verri} and involves the adaptation of tools from algebraic topology into methods that analyze datasets, imagery, and signals. It has found applications in areas as diverse as protein structure \cite{Headd}, gene expression \cite{Deq}, multi-modal fusion \cite{Myers}, time series analysis \cite{Seversky2016}, and road network reconstruction \cite{Ahmed}. In particular, there have been myriad applications of TDA to image analysis, including \cite{Chen2011} and  \cite{VerHoef}.

	More recently, there has been an explosion of connections between TDA and modern machine learning (ML); see \cite{Hensel} for a recent survey. A key element in many TDA/ML pipelines is the computation of a persistence diagram (PD), which is a compact two-dimensional summary of the multi-scale shape of a possibly high-dimensional point cloud (Figure \ref{pers_image_pipeline}). In short, a PD is a multi-set of dots, each of which represents a topological feature (component, hole, void, etc) of a different scale or persistence. PDs can then be vectorized \cite{Adams} and fed into standard ML pipelines, often via PD-bespoke neural network architectures \cite{PersLay}, or they can be used to define shape-based loss functions for use in, for example, segmentation \cite{Hu}.
\begin{figure}
\begin{center}
\begin{tikzpicture}[scale = 0.15,node distance=0.2cm, every node/.style={align=center}]
    \node (n1) {\textbf{Point cloud} \\[2pt] {\centering \includegraphics[width=2.5cm]{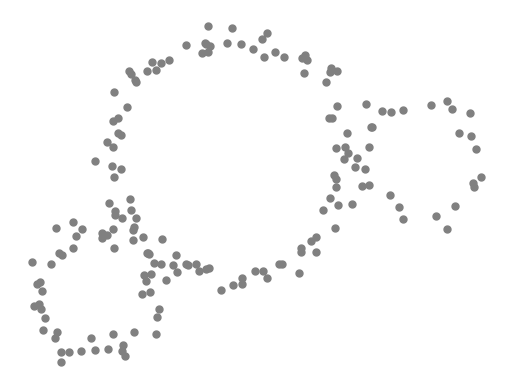}\vspace{2pt}} };
    \node (n2) [below =of n1,xshift = .5cm,yshift = -0.8cm] {\textbf{Cross-Matched}\\\textbf{Weighted Graph} \\ {\centering\includegraphics[scale=.15]{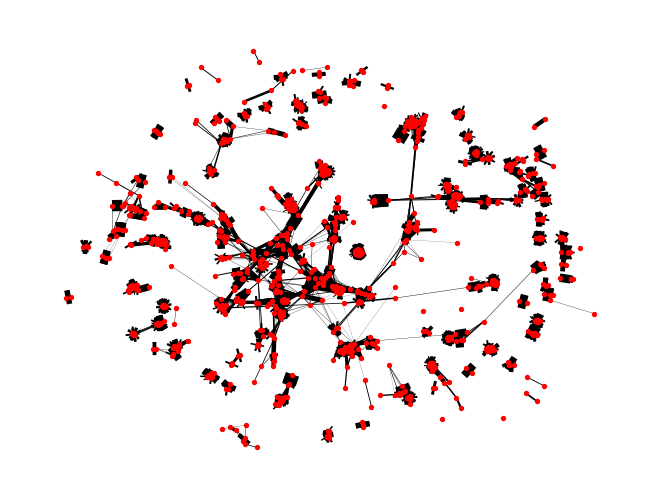}}};
    \node (n3) [below=of n1,xshift = 4.3cm,yshift = 1.5cm] {\textbf{Cross-Matched}\\ \textbf{Prevalence Image} \\[3pt] {\centering\includegraphics[scale=.20]{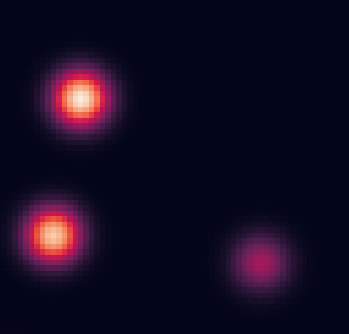}\vspace{3pt}}};

    \draw[->,line width=1pt,solid,bend right] (n1) to node[right,draw = white,xshift = .25cm] {(1) sample and\\ cross-match} (n2);
    \draw[->,line width=1pt,bend right,solid] (n2) to node[below right, draw = white,yshift = -.1cm] {(2) perform graph\\ clustering} (n3);
\end{tikzpicture}
\end{center}
\caption{Cross-Matched Prevalence Image Pipeline}
\label{pipeline}
\end{figure}

\begin{figure}
	\begin{center}
	\begin{subfigure}{0.48\linewidth}
	\centering
    \includegraphics[scale = .25]{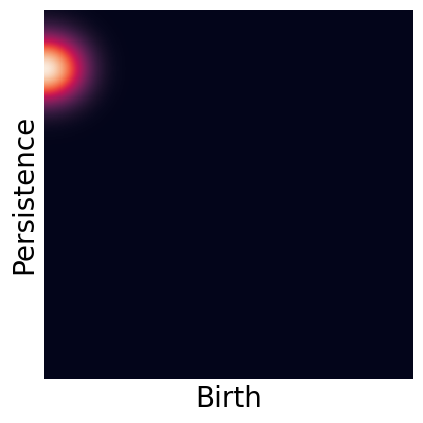}
    \caption{Noiseless}
	\end{subfigure}
    \hfill
    \begin{subfigure}{0.48\linewidth}
    \centering
    \includegraphics[scale = .25]{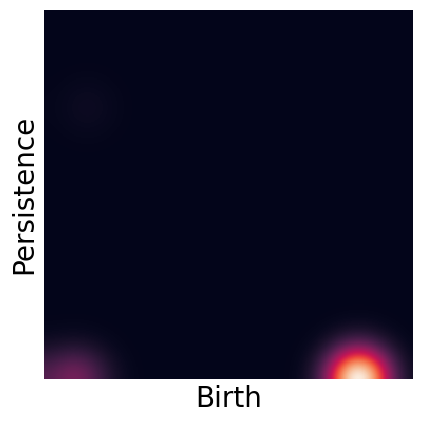}
    \caption{Noise, $\nu = 0.1$}
    \end{subfigure}
    \end{center}
    \caption{Persistence Image of Circle Embedded in High Dimension. The prominent loop indicated by the bright spot on the left is degraded by the addition of noise.}
    \label{pers_image_example}
	\end{figure}

	Typically, PDs in ML are computed for a set of $M$ points in $D$-dimensional ambient (Euclidean) space, $M$ and $D$ being large in both present problems. When $M$ is large, there are significant speed and memory issues for fairly standard linear algebra reasons \cite{Otter}. A more subtle issue arises when points are sampled, with low levels of noise, from a lower-dimensional shape $\mathcal{S}$ within the ambient space; this is often the assumption made by the so-called 'manifold hypothesis' when $\mathcal{S}$ is for example a set of high-resolution images. As Reani and Bobrowski \cite{Reani} observe, even if $\mathcal{S}$ itself has well-resolved shape of significant size, the resulting persistence diagram typically looks like one sampled from a solid ball; in effect, the “curse of dimensionality” manifests here as a destruction of any underlying topological signal (Figure \ref{pers_image_example}).
	
	Reani and Bobrowski propose a solution to this latter problem, namely replacing the idea of persistent homology with that of topological prevalence: in a nutshell, one takes many bootstrap samples of the original point cloud and computes how often the topological shape in the large point cloud is reflected in the bootstrap samples. Their finding, supported with significant experimental evidence by other authors \cite{Easley,GR} is that highly prevalent features correspond to the actual shape of $\mathcal{S}$. However, their solution does not address the first problem; indeed, computation of topological prevalence requires computation of the entire PD of the large point cloud in addition to that of all bootstrap samples.

	In this paper, we propose a method that solves both issues simultaneously. We present a pipeline to approximate topological prevalence of a point cloud $X \subset \mathbb{R}^D$ that does not require a persistent homology computation on the entire space $X$, rather only persistent homology computations on the scale of bootstrap samples of $X$. The pipeline has three primary steps. 1. We perform bootstrap sampling on $X$. Using the machinery of Reani and Bobrowski \cite{Reani} and Garcia-Redondo \etal \cite{GR}, we match topologically prevalent features between $2$ samples at a time. This matching process yields an $N$-partite graph where $N$ is the number of bootstrap samples. 2. We perform a clustering algorithm on the resulting graph. 3. We then compute an image from the output of clustering called a \textit{cross-matched prevalence image (CMPI)} that displays an estimate of the birth and death times of persistent features (see Figure \ref{pipeline}).

%% file: sec/2_prelim.tex
\section{Preliminary}
\label{sec:preliminary}
\subsection{Persistent Homology}
\begin{figure}
\begin{center}
\begin{tikzpicture}[scale = 0.15,node distance=0.2cm, every node/.style={align=center}]
    \node (n1) {\textbf{Point cloud} \\[2pt] {\centering \includegraphics[scale = 0.25]{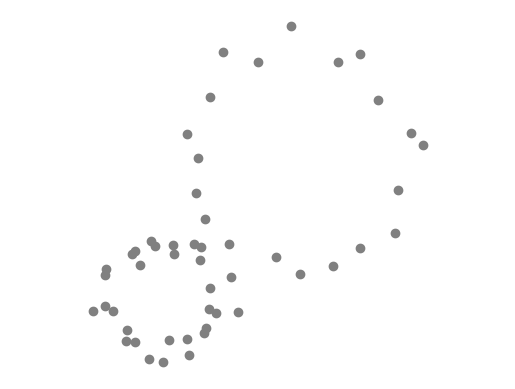}\vspace{2pt}} };
    \node (n2) [below =of n1,xshift = 0cm,yshift = -1.1cm] {\textbf{Persistence Diagram} \\ {\centering\includegraphics[scale=.25]{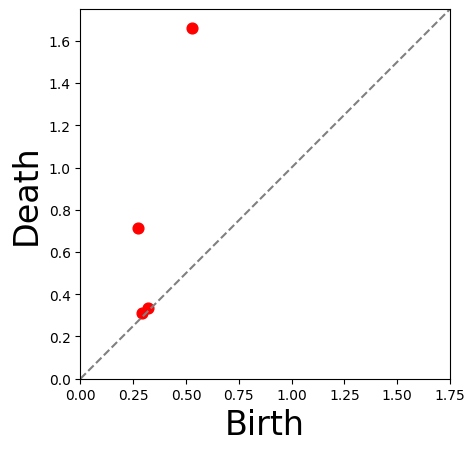}}};
    \node (n3) [below=of n1,xshift = 4.3cm,yshift = 1.5cm] {\textbf{Persistence Image} \\[3pt] {\centering\includegraphics[scale=0.25]{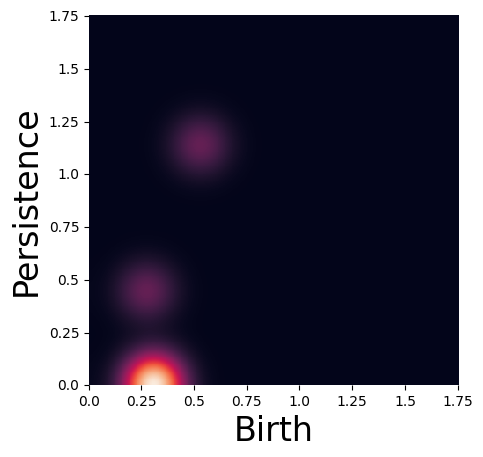}\vspace{3pt}}};

    \draw[->,line width=1pt,solid,bend right] (n1) to node[right,draw = white,xshift = .25cm] {persistent homology\\ computation} (n2);
    \draw[->,line width=1pt,bend right,solid] (n2) to node[below right, draw = white,yshift = -.1cm] {apply $T$, weighting $f$\\ and Gaussian blur} (n3);
\end{tikzpicture}
\end{center}
\caption{Persistence Diagram and Persistence Image}
\label{pers_image_pipeline}
\end{figure}

For a more algebraically rigorous presentation of fundamental tools in persistent homology, we direct interested readers to view \cite{ELZ} and \cite{Verri}. We provide an intuitive overview here. Although the theory can be made more general, here we are studying persistent homology for point clouds $X$ within an ambient space $\mathbb{R}^D$. As the points within $X$ thicken (i.e. replaced by a ball of increasing radius) within the ambient space, topological features (e.g. components, holes, voids, ...) are formed and disappear. The \textit{persistence diagram} $\dgm(X)$ summarizes this process containing a dot $(b,d)$ for each feature that appears at radius $b$ and disappears at radius $d$. We often refer to $b$ and $d$ as \textit{birth} and \textit{death} times respectively. See Figure \ref{pers_image_pipeline}.

We call the dots $(b_i,d_i)$ contained in a persistence diagram the \textit{bars} or \textit{intervals} of $X$ and the multiset $\{(b_i,d_i)\}$ of all bars in a given persistence diagram is the \textit{persistence barcode} or \textit{barcode} $\barcode(X)$.

\begin{remark}
Throughout this paper, in all examples and experimentation we compute degree $1$ persistent homology which encodes the information of holes. All algorithms discussed can be extended to other degrees, such as degree $0$ which encodes information of components and degree $2$ which encodes information of voids.
\end{remark}

\subsubsection{Persistence Image}
\label{pers_im_sec}
The stable vectorization of a barcode is an important area of study in persistent homology. Such schemes, for example, make persistent homology a viable candidate as a feature-extractor in the pre-processing step of machine learning pieplines. The \textit{persistence image} is one such stable vectorization \cite{Adams}.

\begin{definition}
Given a weighting function $f \colon \mathbb{R}^2 \to \mathbb{R}_{\geq 0}$, the \textit{persistence image} associated with a barcode $\{(b_i,d_i)\}_{i \in I}$ is a discretization of the $2$-dimensional surface
\begin{equation}
\rho(x,y) = \sum_{i \in I} f(\alpha_i) g_{\alpha_i}(x,y)
\end{equation}
where $T\colon \mathbb{R}^2 \to \mathbb{R}^2$ is the linear transformation $(x,y) \mapsto (x,y-x)$, $\alpha_i \in \mathbb{R}^2$ is $T(b_i,d_i)$ and $g_{\mu} \colon \mathbb{R}^2 \to \mathbb{R}$ is a normalized symmetric Gaussian with mean $\mu$ and variance $\sigma^2$. 
\label{def_pers_image}
\end{definition}

\begin{remark}
Because the horizontal axis of a persistence image is informed by the birth of intervals and the vertical axis by persistence of intervals, we say that the persistence image is on the \textit{birth-persistence plane}. Note a persistence diagram is on the \textit{birth-death plane}.
\end{remark}

Under reasonable assumptions on $f$ \cite{Adams}, the persistence image is a stable representation of a persistence diagram, making it a viable candidate as a feature-extractor in the pre-processing step of machine learning pipelines. 

\subsection{Image Persistence and Interval Matching}
\label{prev_preliminary}
Given a point cloud $C$ and a subset $A \subset C$, the \textit{image persistence} of $A$ in $C$ contains the data of topological features also present in $C$ as both point clouds thicken.
Given two subsets $A,B \subset C$, image persistence provides the machinery to quantify the similarity of bars of $A$ and bars of $B$ through $C$. Image persistent homology, like its traditional counterpart, outputs a multiset of bars. Let $\barcode_C(A)$ denote the \textit{image barcode} for $A$ in $C$.  This is made rigorous by the concept of interval matching and matching affinity defined below. The interested reader is encouraged to review \cite{Bauer} and \cite{GR} for further details.

\begin{definition}[Interval Matching]
Let $A,B \subset C$. We say the interval $\alpha \in \barcode(A)$ is \textit{matched} with $\beta \in \barcode(B)$ if there exists $\tilde{\alpha} \in \barcode_C(A)$ and $\tilde{\beta} \in \barcode_C(B)$ such that $\text{birth}(\alpha) = \text{birth}(\tilde{\alpha})$, $\text{birth}(\beta) = \text{birth}(\tilde{\beta})$, and $\text{death}(\tilde{\alpha}) = \text{death}(\tilde{\beta})$.
\label{interval_matching}
\end{definition}

\begin{definition}[Matching Affinity]
\newcommand{\J}{\mathcal{J}}
Let $A,B \subset C$. For $\alpha \in \barcode(A)$ matched with $\beta \in \barcode(B)$ with $\tilde{\alpha}$ and $\tilde{\beta}$ as in the previous definition, the associated \textit{matching affinity} is the quantity
\begin{equation}
\affinity(\alpha,\beta) = \J(\alpha,\beta) \cdot \J({\tilde \alpha}, \alpha) \cdot \J({\tilde\beta},\beta)
\label{eq:aff}
\end{equation}
where $\J$ is the Jaccard index defined by $\J(I,J) := |I \cap J| / |I \cup J|$.
\end{definition}

The notion of interval matching is made particularly useful when we take a bootstrap sample $\{X_1, \dots, X_N\}$ of a point cloud $X$. Image persistence then enables the identification of topological features of $X$ present across many bootstrap samples $X_i$ by setting $A = X_i$, $B = C = X$ in the above definitions.

\begin{remark}
Note for $B = C = X$, $\barcode_C(B) = \barcode(B)$, thus $\tilde{\beta} = \beta$, simplifying the statements in Definition \ref{interval_matching}.
Further, each bar $\alpha \in \barcode(X)$ has a unique match in $\barcode(X_i)$ which we denote as $\beta_i(\alpha)$ \cite{GR}. The latter remark allows for the following construction.
\end{remark}

\begin{definition}[Prevalence Score]
Given a sample $\{X_1, \dots, X_N\}$ of $X$, and $\alpha \in \barcode(X)$, the \textit{prevalence score} of $\alpha$ is the quantity
\begin{equation}
\text{prev}(\beta) = \frac{1}{N} \sum_{i = 1}^N \affinity(\alpha,\beta_i(\alpha)).
\end{equation}
\label{prev_score}
\end{definition}

\begin{remark}
A bar $\alpha \in \barcode(X)$ need not have a match in all (or any) samples. In such cases in the above definition, $\beta_i(\alpha)$ may be ill-defined, thus by convention we set any associated summands to $0$.
\end{remark}
Bars $\alpha \in \barcode(X)$ receiving relatively high prevalence scores are said to be \textit{prevalent}. A primary theme to be observed in the discussion that follows is that the notion of prevalence is more robust to noise in high dimensional settings then is persistence.

%% file: sec/3_cm.tex
\section{Cross-Matching}
\label{sec:cm}
We now work in the following context.  Fix a finite point cloud $X \subset \mathbb{R}^D$ and a set of bootstrap samples $\mathcal{S} = \{X_j \subset X \,\mid\, 1\leq j \leq N,\,|X_j| = s\}$. By a \textit{reference bar}, we mean an element of $\barcode(X)$ and by \textit{sample bar} we mean an element of $\barcode(X_i)$ for some $j = 1, \dots, N$.

The primary contribution of this paper is an algorithm which approximates the prevalence information of a point cloud $X \subset \mathbb{R}^D$ via bootstrap sampling which does not require the computation of $\barcode(X)$ thereby bypassing computational difficulty when $X$ is large in size. The main export of this algorithm is the \textit{cross-matched prevalence image} (see section \ref{cm_image}), an image on the birth-persistence plane that displays the prevalence of cross-matched features (Definition \ref{CM_graph_def}) on the scale of sample bars. To assess the accuracy of this image, we introduce an image called the \textit{prevalence image} in section \ref{prev_image} that summarizes the information of reference bars, their prevalence scores, and their matches across samples. Note as the prevalence image is dependent on the information of traditional interval matching, it requires the computation of $\barcode(X)$ and thus for some point clouds $X$ may not be feasibly computed (see Section \ref{complexity}). The prevalence image is defined in such a way that its features are on the scale of sample bars, not on the scale of reference bars. This choice is necessary for the prevalence image and cross-matched prevalence image to be comparable.

\subsection{Prevalence Image}\label{prev_image}
Inspired by the persistence image (Section \ref{pers_im_sec}), we introduce a new image called the \textit{prevalence image} which given a choice of bootstrap sample encodes the prevalence of topological features of a space in the birth-persistence plane. 

A prevalence image is constructed as follows. We wish to define a vector $\mu_i \in \mathbb{R}^2$ for each reference bar $I_i$ to be plotted on the birth-persistence plane that is on the scale of sample bars rather than on the scale of reference bars. Let $I_i^{(j)} \in \barcode(X_j)$ denote a match with $I_i$ with affinity $\alpha_i^{(j)}$ for $j = 1, \dots , N$. Corresponding to each reference bar $I_i$, define a discrete probability space $(\Omega_i, 2^{\Omega_i}, \mathbb{P}_i)$ where $\Omega_i = \{I_i^{(j)}\}$ and $\mathbb{P}_i(I_i^{(j)}) = \alpha_i^{(j)} / \left(\sum_{k = 1}^N \alpha_i^{(k)}\right)$. Let $T$ and $g_\mu$ be as in Definition \ref{def_pers_image}. Then for each reference bar $I_i$ we define a vector $\mu_i \in \mathbb{R}^2$ to be the expectation $\mathbb{E}[T] := \int_{\Omega_i} T d\mathbb{P}_i$.
\begin{definition}[Prevalence Image]
For $X \subset \mathbb{R}^d$ a finite point cloud and a set of bootstrap samples $\mathcal{S}$, the \textit{prevalence image} is a discretization of the surface
\begin{equation}
\rho^{(X,S)}_\prev(x,y) = \sum_{ I_i \in \barcode(X)} \prev(I_i) \,g_{\mu_i}(x,y).
\end{equation}
\label{prev_surf}
\end{definition}
\begin{remark}
Recall the reference bar $I_i$ need not have a match in each $\barcode(X_j)$ and in such an instance $\alpha_i^{(j)} = 0$ by definition. Consequentially, in the event $I_i$ has no matches in any sample with nonzero matching affinity, $\prev(I_i) = 0$. Thus, we adopt the convention that any summand in Definition \ref{prev_surf} with $\prev(I_i) = 0$ is disregarded as $\mu_i$ is ill-defined in such cases.
\end{remark}

\subsection{Cross-Matching Intervals and Graphs}
Our proposed approximation for the prevalence image, the aforementioned cross-matched prevalence image, depends entirely on a matching scheme to be discussed now in detail. The conventional approach to identifying topologically prevalent intervals as pursued in the recent work \cite{Reani, GR, Easley} matches reference bars with sample bars through image persistence, as described in section \ref{prev_preliminary}. Such a matching scheme requires knowledge of (and therefore computation of) $\barcode(X)$ (see Definition \ref{prev_score}).

Here we propose an alternative matching scheme, called \textit{cross-matching}, and the related \textit{cross-matched weighted graph}. The name \textit{cross-matching} itself was inspired by a method  of the same name given in the publicly available python repository of Garcia-Redondo \etal \cite{GR}. However, as far as we are aware the concept itself has not been studied to date as an alternative to conventional matching in the context of topological prevalence.
\begin{definition}[Cross-Matched Intervals and Cross-Matching]
Let $X,Y \subset \mathbb{R}^d$ be finite point clouds. We say that $\alpha \in \barcode(X)$ and $\beta \in \barcode(Y)$ are \textit{cross-matched} if they are matched through their union $X \cup Y$. That is, they are matched according to Definition 2.2 with $A = X, B = Y$ and $C = X \cup Y$. We write $\alpha \sim_{\text{CM}} \beta$ for $\alpha$ and $\beta$ cross-matched. For a set of samples $\mathcal{S} = \{X_1, \dots, X_N\}$, by \textit{cross-matching}, we mean identifying all cross-matched intervals of $\barcode(X_i)$ and $\barcode(X_j)$ for all distinct $i$ and $j$.
\end{definition}
\begin{definition}[Cross-Matched Weighted Graph]
For a set $\mathcal{S} = \{X_1, \dots, X_N\}$ of samples, let $V_\mathcal{S}$ denote the set of all sample bars $I \in \barcode(X_i)$ for some $i$ that are cross-matched with some sample bar $J \in \barcode(X_j)$ for $j \neq i$. Then let $E_\mathcal{S} \colon V_\mathcal{S} \times V_\mathcal{S} \to [0,1]$ where $E_\mathcal{S}(I,J) = \affinity(I,J)$ if $I$ and $J$ are cross-matched, else it is $0$. Then let $G_\mathcal{S}$ be the weighted graph with vertex set $V_\mathcal{S}$ and edge weights given by $E_\mathcal{S}$. We call $G_\mathcal{S}$ the \textit{cross-matched weighted graph} of $\mathcal{S}$.
\label{CM_graph_def}
\end{definition}
\begin{remark}\label{partite_remark}
Because in the process of cross-matching, comparisons are made only between distinct sampled point clouds $X_i$ and $X_j$, for $\mathcal{S}$ of size $N$, the cross-matched weighted graph $G_\mathcal{S}$ is $N$-partite.
\end{remark}
\subsection{Graph Clustering}

We wish to employ clustering schemes on the weighted $N$-partite graph $G_{\mathcal{S}}$ to extract prevalence information of the original point cloud $X$. In pursuit of this objective, we first distinguish two factors that impact the prevalence score of a reference bar $I$ (see Definition \ref{prev_score}): (1) the number of samples $X_j$ with a sample bar matched to $I$, or the \textit{frequency} of matches across the set $\mathcal{S}$, and (2) the value of matching affinities of said matches, or the \textit{quality} of matches across the set $\mathcal{S}
$. One might wish to then adopt a scheme to extract information from $G_{\mathcal{S}}$ that is informed both by the \textit{frequency} and the \textit{quality} of matches with reference bars. However, no such scheme is feasible as by its very construction no explicit reference bar information is present in $G_{\mathcal{S}}$. So, instead we seek a clustering scheme that is informed by the frequency and quality of \textit{cross-matches} among sample bars. Thus, the clustering scheme adopted ought to consider not only the edge set (frequency) but also the edge weights (quality) of $G_{\mathcal{S}}$ when assigning vertices to clusters.

\begin{definition}[Cluster Inter-degree and Intra-degree]
Let $G = (V,E)$ be a weighted graph, and let $\mathcal{C}$ be a clustering scheme of $G$. For $C \in \mathcal{C}$ and $v \in C$, let $d_C(v)$ denote the sum of all edge weights between $v$ and any other vertex $w \in C$. Let $\overline{d_C}(v)$ denote the sum of all edge weights between $v$ and any other vertex $w \in V \cap C^C$. Then $d_C(v)$ and $\overline{d_C}(v)$ are, respectively, the \textit{cluster inter-degree} and \textit{cluster intra-degree} of $v$ in $C$. Define the \textit{degree of a cluster}, $d_C$, to be half the sum of $d_C(v)$ over all $v \in C$.
\end{definition}
\begin{remark}
The degree of a cluster $C$ can equivalently be defined to be the sum of all edge weights between vertices both contained within $C$ with each edge being counted only once.
\end{remark}

We implement a greedy clustering algorithm informed by both the edge set (frequency) and edge weights (quality) that is designed according to the following objectives (Section \ref{alg_gc}). (1) The quantity $\sum\limits_{C \in \mathcal{C}} d_C$ is large. (2) The quantity $\sum\limits_{C \in \mathcal{C}} \sum\limits_{v \in C} d_C(v)$ is small.
\subsection{Cross-Matched Prevalence Image}
\label{cm_image}
We are now equipped to define the cross-matched prevalence surface and image. Let $\mathcal{C}$ be a clustering scheme on $G_{\mathcal{S}}$. For $v$ a vertex of $G_{\mathcal{S}}$, let $I_v$ denote the associated sample bar. Now, for each cluster $C \in \mathcal{C}$, define a discrete probability space $(\Omega_C,2^{\Omega_C},\mathbb{P}_C)$ where $\Omega_C = \{ I_v \,\mid\, v \in C\}$ and $\mathbb{P}_C(I_v) = d_C(v) / 2d_C$. For each cluster $C \in \mathcal{C}$, we assign a representative $\mu_C \in \mathbb{R}^2$ defined by $\mu_C := \mathbb{E}[T] =  \int_{\Omega_C}T\, d \mathbb{P}_C$. 

\begin{definition}[Cross-Matched Prevalence Image]
Let $X \subset \mathbb{R}^D$ be a finite point cloud and $\mathcal{S}$ a set of bootstrap samples of $X$. Fix a clustering scheme $\mathcal{C}$ on $G_{\mathcal{S}}$. The \textit{cross-matched prevalence image} is a discretization of the surface
\begin{equation}
\rho_{\text{CM}}^{(X,\mathcal{S},\mathcal{C})}(x,y) = \sum_{C \in \mathcal{C}} d_C\, g_{\nu_C}(x,y).
\end{equation}
\end{definition}

\begin{remark}
Unless otherwise stated, anytime we refer to a cross-matched or reference prevalence image, we are referring to the image after $L^1$ normalization, that is dividing the value at each pixel by the sum of the absolute value of all pixels. 
\label{norm_remark}
\end{remark}
\begin{remark}
By construction of the prevalence and cross-matched prevalence images, regions with high value are regions associated with highly persistent intervals. Plots of images in this paper are colored ranging in spectrum from yellow (highest value) to purple (lowest value). By a feature in an image, we mean a well-isolated bright spot. Brightness is correlated with persistence directly. When we say a feature is dim and dimming with $s$, we mean the prevalence value is relatively low and the prevalence value is decreasing as $s$ increases (respectively).
\end{remark}

%% file: sec/4_implementation.tex
\section{Implementation}
\label{algorithms}

\subsection{Cross-Matching}
\label{alg_cm}
We now detail an algorithm which computes the Cross-Matched Weighted Graph (Definition \ref{CM_graph_def}) for a point cloud $X$ given a bootstrap sample $\mathcal{S}$.
\begin{algorithm}

\caption{Compute Cross-Matched Weighted Graph on Point Clouds}
\label{alg_1}
\begin{algorithmic}[1]
\Require $X$, point cloud of size $M$ in $\mathbb{R}^D$ as a $M \times d$ array
\Require $\mathcal{S} = \{X_1, X_2, \dots, X_N\}$, a set of samples of $X_i \subset X$ with $|X_i| = s \in \mathbb{N}$ for all $i$
\Procedure{CMGraph}{$X$,$\mathcal{S}$}
\State Initialize empty graph $G$
\For{$X_i \in \mathcal{S}$}
    \State Compute $\barcode(X_i)$
\EndFor
\For{$X_i \in \mathcal{S}$}
	\For{$X_j \in \mathcal{S}$ with $j > i$}
		\State $E_{i,j} \gets \{ (\alpha,\beta) \in \barcode(X_i) \times \barcode(X_j) \mid \alpha \sim_{CM} \beta\}$
		\State $V_{i,j} \gets \{\text{bars of } X_i \text{ matched}\} \cup \{\text{bars of } X_j \text{ matched}\}$
	\EndFor
\EndFor
\State $E \gets \bigcup_{j > i} E_{i,j}$
\State $V \gets \bigcup_{j > i} V_{i,j}$
\State Set the vertex set of $G$ to $V$
\For{$(\alpha,\beta)$ in $E$}
	\State Add edge between $\alpha$ and $\beta$ with weight $\affinity(\alpha,\beta)$ to $G$
\EndFor
\EndProcedure
\end{algorithmic}

\end{algorithm}
\begin{remark}
In the above algorithm, the contents of the for-loop spanning lines 4-5 can be ran in parallel using up to $N$ workers, and the contents of the nested for-loop spanning lines 6-10 can be ran in parallel using up to $\binom{N}{2}$.
\end{remark}
\begin{remark}
We use routines in the repository of Garcia-Redondo \etal \cite{GR} for operations in Lines $4$ and $8$.
\end{remark}
\subsubsection{Computational Complexity}
\label{complexity}
The computational complexity of Algorithm \ref{alg_1} is on the order of the calculations in lines $4$ and $8$. Line $4$ consists of persistent homology computations while line $8$ consists of image persistent homology computations. To compute persistence homology, one uses the Vietoris-Rips (VR) complex. In a nutshell, the computation of degree $i$ persistent homology (both standard and image) requires the reduction of matrices that are of size corresponding to the number of simplices up to degree $i + 1$ in the VR complex associated to the point cloud of interest. In the worst case, the VR complex for a point cloud of size $M$ has $\binom{M}{k+1}$ $k$-simplices. For degree $1$ persistent homology, we must consider $k$-simplices for $k = 0,1,2$. So for a point cloud $X_i$ of size $s$, we must reduce a matrix of size on the order of $\mathcal{O}(\binom{s}{1})+ \mathcal{O}(\binom{s}{2})+\mathcal{O}(\binom{s}{3}) \sim \mathcal{O}(\binom{s}{3})$ which is complexity $\mathcal{O}(\binom{s}{3}^3) \sim \mathcal{O}(s^9)$. Then, lines 3-5 have complexity $\mathcal{O}(Ns^9)$. The matching algorithm from the repository of Garcia-Redondo \etal \cite{GR} used in line 8 requires the computation of image persistent homology of clouds $X_i$ and $X_j$ within $X_i \cup X_j$. A degree $i$ image persistence computation of a point cloud $B \subset A$ requires the reduction of two matrices, both of size on the order of $\mathcal{O}(\binom{|A|}{i+1})$ where $|A|$ is the size of the point cloud $A$. For further technical details, see \cite{Bauer}. Then in our setting, for degree $1$ and given $A = X_i \cup X_j$ is size $2s$, line $8$ is complexity $\mathcal{O}(\binom{2s}{3}^3)$. Line 8 is run $N^2$ times, thus in total the complexity for the for-loop spanning lines 6-10 is $\mathcal{O}(N^2 s^9)$. Then the overall complexity is $\mathcal{O}(N^2 s^9) + \mathcal{O}(N^2 s^9) \sim \mathcal{O}(N^2 s^9).$ 

Note, as the conventional prevalence scheme requires the computation of $\barcode(X)$, its complexity is $\mathcal{O}(M^9)$ where $|X| = M$. In practice, $s << M$. Assuming access to parallel workers, the cross-matched graph can be computed with complexity as small as $\mathcal{O}(s^9)$. 
\subsection{Graph Clustering}
\label{alg_gc}
Given a weighted $N$-partite graph $G = (V,E)$, we employ an intuitive greedy clustering algorithm to form clusters of high cluster degree. The algorithm begins by forming a queue $Q$ from $V$ by sorting vertices by weighted degree in descending order. The algorithm proceeds through the $Q$ forming clusters of the current first in line vertex and its neighbors. After making this initial choice of clustering, we then move through the set of all neighbors of this newly formed cluster $C$ (also by descending degree), and replace any vertex $v \in C$ with a neighbor of the cluster $n$ if the resulting cluster degree increases by doing so. See Algorithm \ref{alg_2} for a formal statement of the algorithm.
\newcommand{\cdeg}[1]{\text{deg}(#1)}
\begin{algorithm}
\caption{Graph Clustering}
\label{alg_2}
\begin{algorithmic}[1]
\Require $G$, weighted $N$-partite graph with vertex set $V$ and edge set $E$
\Procedure{HDCluster}{$G$}
\State $\mathcal{C} \gets \text{empty list}$
\State $Q \gets \text{sort } V \text{ by weighted degree (descending)}$
\While{$|Q| > 1$}
\State $x \gets \text{pop}(Q)$ 
\State $C \gets \{x\} \cup \{\text{neighbors of } x\}$
\State $Q_n \gets \text{NeighborQueue}(\{\text{neighbors of } x\},Q)$
\While{$Q_n$ is non-empty}
\State $n \gets \text{pop}(Q_n)$
\If{$\exists\, v \in C$ with $v$ in partition of $n$}
\If{$\cdeg{C}< \cdeg{C \backslash \{v\} \cup \{n\}}$}
\State Replace $v$ with $n$ in $C$
\EndIf
\EndIf
\EndWhile
\For{$y \in C$}
\State Remove $y$ from $Q$
\EndFor
\State $\mathcal{C} \gets \mathcal{C} + \{C\}$
\EndWhile
\State \Return $\mathcal{C}$
\EndProcedure
\Require $\mathcal{N}$, subset of vertex set $V$
\Require $Q$, queue containing vertices of $V$
\Procedure{NeighborQueue}{$\mathcal{N}$,$Q$}
\State $Q_{n} \gets \text{empty queue}$
\For{$n \in N$}
\For{$y \in \{\text{neighbors of } n \text{ in }  Q\}$}
\State $Q_n \gets Q_n + \{y\}$
\EndFor
\EndFor
\State \Return $Q_n \text{ sorted by descending weighted degree}$
\EndProcedure
\end{algorithmic}
\end{algorithm}

%% file: sec/5_param.tex
\section{Selection of Sample Parameter $s$}
\label{sec:param}
Given a point cloud $X \in \mathbb{R}^D$ and a bootstrap sample $\mathcal{S} = \{X_1, \dots, X_N\}$, one can complete the cross-matched prevalence image using Algorithm \ref{alg_1} to obtain the associated cross-matched weighted graph $G_\mathcal{S}$ and the clustering algorithm (Section \ref{alg_gc}) to identify high degree clusters. To proceed in this manner, one must select two parameters, the size $s$ of each $X_i$, and the total number $N$ of samples $X_i$. Particularly in noisy settings, $s$ has a substantial effect on the resulting image. It is therefore necessary to select such a parameter with care. In this section, we propose an early-stopping algorithm to select $s$ that can be used for both cross-matched and reference images.

\begin{figure}
\subcaptionsetup{labelformat=empty}
    \centering
    \begin{minipage}[t]{0.60\linewidth}
        \centering
        \vfill
        \includegraphics[width=0.8\linewidth]{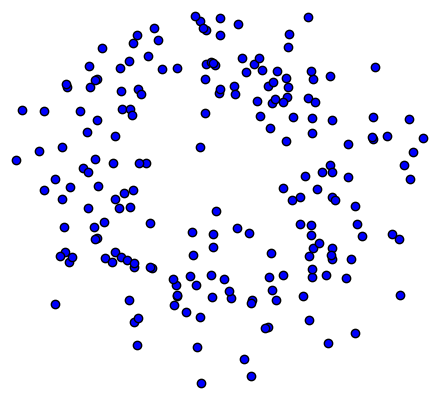} 
        \subcaption{Noisy circle}
        \vfill
    \end{minipage}%
    \hspace{1cm}
    \begin{minipage}[t]{0.15\linewidth}
        \centering
        \vfill
        \includegraphics[width=0.6\linewidth]{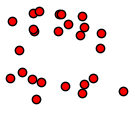} 
        \subcaption{$s = 25$}
        \vspace{1em}
		\vfill
        \includegraphics[width=0.6\linewidth]{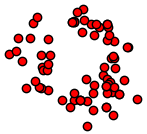} 
        \subcaption{$s = 75$}
        \vspace{1em}
		\vfill
        \includegraphics[width=0.6\linewidth]{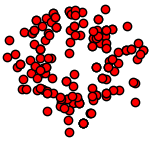} 
        \subcaption{$s = 175$}
        \vspace{1em}
        \vfill
    \end{minipage}
    \captionsetup{labelformat=default}
    \caption{Noisy circle of size $200$ and samples of varying size $s$}
    \label{param_selection_cloud}
\end{figure}

\begin{figure}
    \centering
    \includegraphics[width = \linewidth]{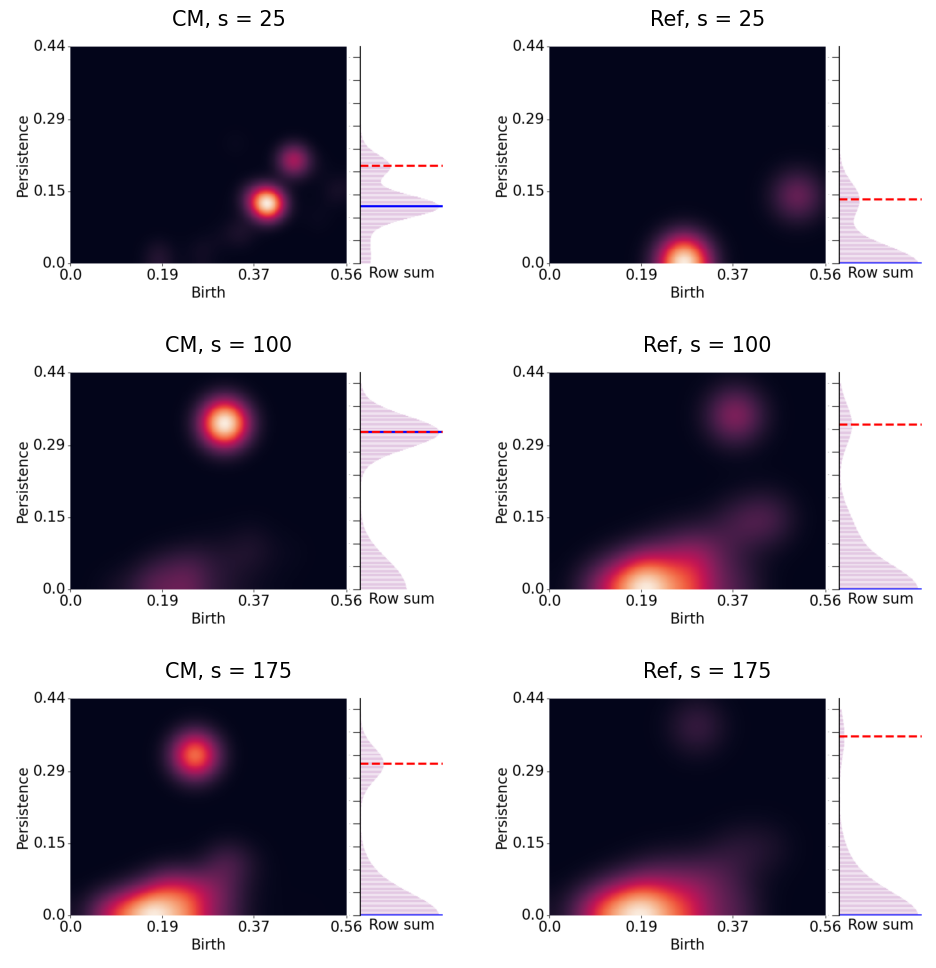}
    \caption{Cross-Matched and Reference Images for noisy circle across varying sample sizes $s$}
    \label{s_param_selection}
\end{figure}

In Figure \ref{s_param_selection}, we display the cross-matched and reference images for a noisy circle consisting of $200$ points in $\mathbb{R}^2$ (see Figure \ref{param_selection_cloud}) for increasing value of $s$. Given that a circle has $1$ prominent loop and these images are approximating degree $1$ homology, a \textit{good} image in this context would be one with a single well-isolated feature with high persistence and early birth. It is observable in general that for $s$ values too low and for $s$ values too high, the corresponding images are poor. We observe that initially as $s$ increases, features tend to increase in persistence and in the reference case and the most persistent features increase in prevalence (compare $s = 25$ with $s = 100$ for CM). As $s$ nears $200$ (the size of the entire point cloud), the persistence of features decreases and in particular the prevalence of the most persistent features decreases (compare $s = 100$ and $s = 175$). This is intuitively to be expected. For a noisy cloud, as $s$ increases, it becomes more likely that samples contain outliers that decrease the radius of loops during thickening or remove them completely  (compare $s = 25$ to $s = 175$ in Figure \ref{param_selection_cloud}).

\begin{figure}[h]
    \centering
    \includegraphics[width = 0.75\linewidth]{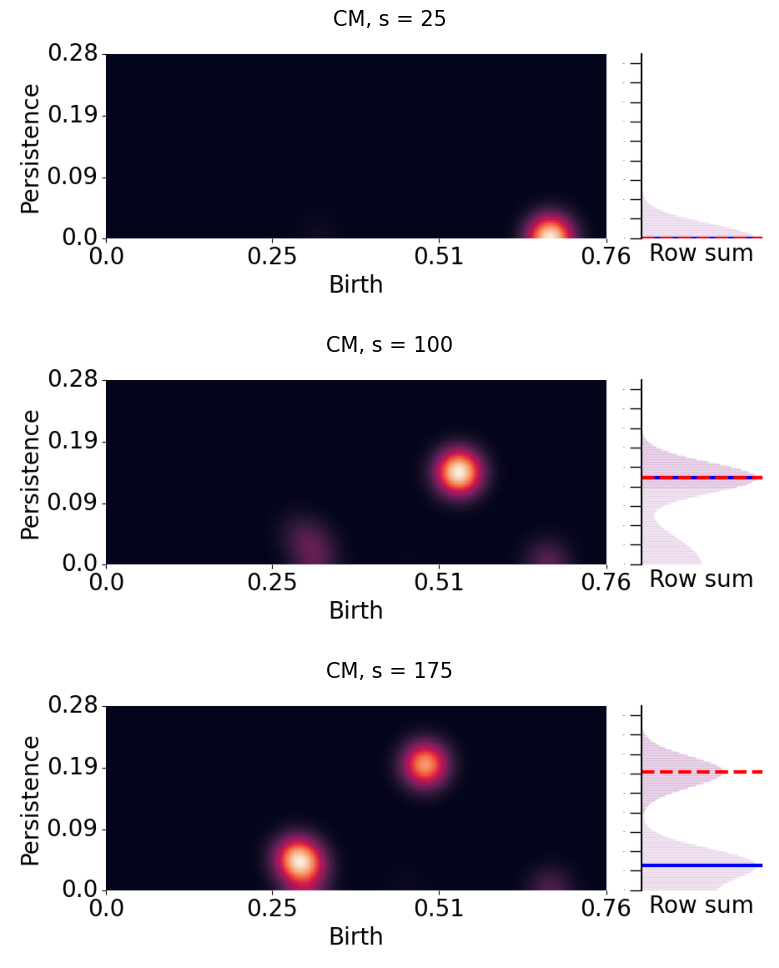}
    \caption{Cross-Matched Images for $T$ for bootstrap samples $S_s$}
    \label{torus_clean_cm_images}
\end{figure}

In Figure \ref{torus_clean_cm_images}, we display cross-matched images for a torus embedded in high dimensional space (see Section \ref{sec:exp} for further details). In this setting, because the degree $1$ homology of a torus is $2$-dimensional, a \textit{good} image would display $2$ prominent well-separated features of similar prevalence. As in the case of the circle, for low $s$ values, images are poor. As $s$ increases, for a time, features become more persistent. From $s = 25$ to $s = 100$, a new feature appears (corresponding to the second generator of degree $1$ homology for a torus). This feature dims the more persistent feature slightly due to the image being normalized (see Remark \ref{norm_remark}) and as $s$ increases further, the two features near the same prevalence.

We propose the following stopping conditions (and stop increasing $s$ when either is first satisfied). (1) The most prevalent feature decreases in persistence. (2) The most persistent feature decreases in prevalence \textbf{while} it is not the most prevalent.

Criteria $1$ prevents sample size from becoming so large that outliers begin to degrade shape (as was the case for the noisy circle for $s$ beyond $75$). Criteria $2$ is related to phenomena that arise when new features are formed, particularly the decrease in prevalence of highly persistent features. In the circle case, these new features are due to noise having greater impact with larger $s$. This is in stark contrast to the torus case, where the prevalence of the most persistent feature decreased due to rising prevalence of a legitimate topological feature. Criteria $2$ balances allowing for the "development" of lower persistent features while preventing the loss of prevalence of highly persistent features to the prevalence of features due to noise.

To implement these stopping criteria, given an image $I$, we consider the row sums of $I$. We define the persistence of the most persistent feature of $I$ to be the height of the highest local maximum of row sums and its prevalence to be the value at said local maximum. We define the most prevalent feature to be the global maximum of row sums. In Figures \ref{s_param_selection} and \ref{torus_clean_cm_images}, we display row sums to the left of each image and mark the most prevalent feature and most persistent feature with blue and red lines respectively.

%% file: sec/6_exp.tex
\section{Experiment: Torus in High Dimensions}
\label{sec:exp}
\newcommand{\imageplotwidth}{0.90\linewidth}
We now compute prevalence and cross-matched prevalence images for a densely sampled $2$-torus $T$ embedded in $\mathbb{R}^{64 \times 64}$. We utilize the publicly available \verb|ellipse| python package that embeds tori, circles, and lines in high dimensional space using parameterizations of ellipses \cite{ellipse}.

Throughout this experiment, $|T| = 484$ and we fix bootstrap samples $\mathcal{S}_{s}$ of $T$ for $s = 25$ through $s = 225$ with a step-size $\Delta s = 25$ such that $|\mathcal{S}_{s}| = 50$ (i.e. $N = 50$) and each $T_i \in \mathcal{S}_{s}$ is sampled with replacement from $T$ with $|T_i| = s$. $T$ is normalized such that each $x \in T$ has $|x| \leq 1$. For all images computed, we fix $\sigma = 0.05$ and compute images of resolution $500 \times 500$. 

\subsection{Accuracy of Cross-Matched Prevalence}

\begin{figure}
\begin{subfigure}{\linewidth}
\centering
\includegraphics[width=\imageplotwidth]{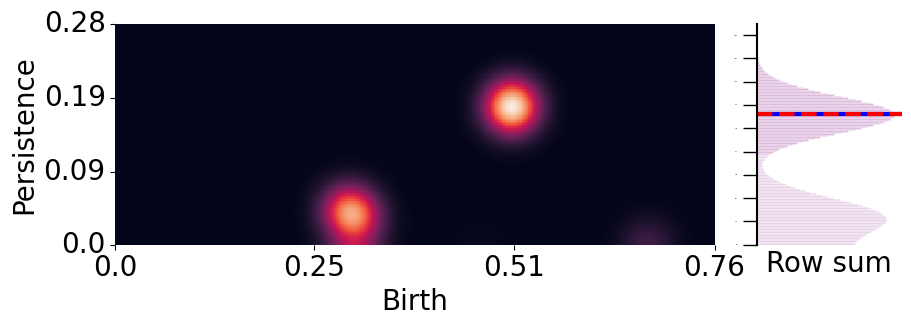}
\caption{CM $(s = 150)$}
\end{subfigure}
\begin{subfigure}{\linewidth}
\centering
\includegraphics[width=\imageplotwidth]{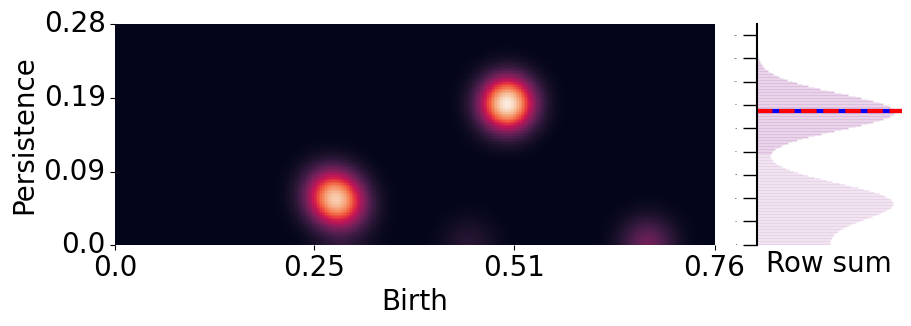}
\caption{Reference $(s = 150)$}
\end{subfigure}
\caption{Persistence Images (Torus, noiseless)}
\label{torus_result_0}
\end{figure}

Given in Figure \ref{torus_result_0} is the output of parameter selection for $s$ (Section \ref{sec:param}) to $T$ for both cross-matched and reference prevalence images for bootstrap sample $\mathcal{S}_{s}$. Both images appropriately depict $2$ prevalent intervals consistent with the degree $1$ homology of a torus. The images agree on approximate birth time and persistence of these $2$ intervals, about $(0.25,0.05)$ and $(0.50,0.20)$ respectively. They also are consistent in identifying the most persistent feature as the most prevalent.

\begin{figure}
\centering
\begin{subfigure}{\linewidth}
\centering
\includegraphics[width=\imageplotwidth]{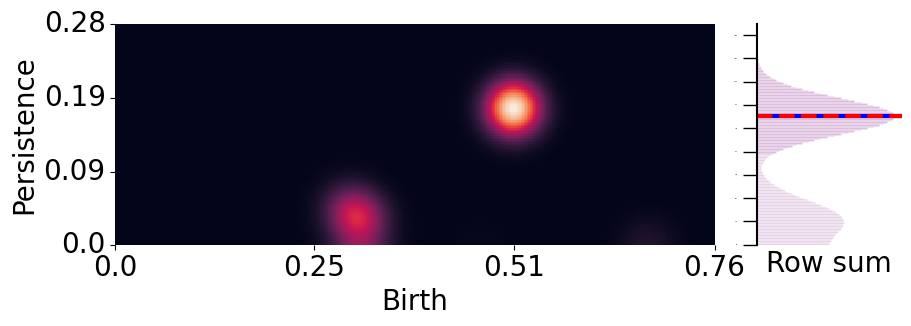}
\caption{CM $(s = 175)$}
\end{subfigure}
\begin{subfigure}{\linewidth}
\centering
\includegraphics[width=\imageplotwidth]{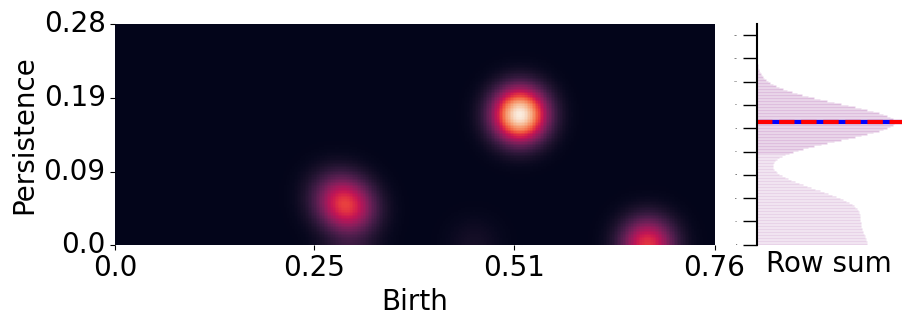}
\caption{Reference $(s = 125)$}
\end{subfigure}
\caption{Persistence Images (Torus, $\nu = 0.1$)}
\label{torus_result_0.1}
\end{figure}

\begin{figure}
\centering
\begin{subfigure}{\linewidth}
\centering
\includegraphics[width=\imageplotwidth]{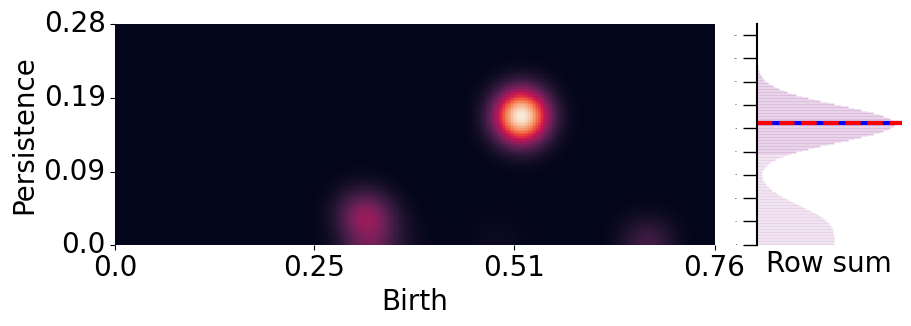}
\caption{CM $(s = 200)$}
\end{subfigure}
\begin{subfigure}{\linewidth}
\centering
\includegraphics[width=\imageplotwidth]{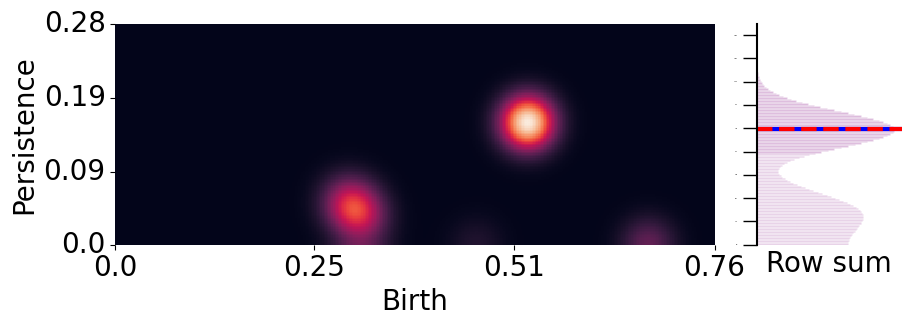}
\caption{Reference $(s = 150)$}
\end{subfigure}
\caption{Persistence Images (Torus, $\nu = 0.3$)}
\label{torus_result_0.3}
\end{figure}

We also introduce Gaussian noise to $T$ with variances $\nu \in \{0.1,0.2,0.3,0.4\}$ to study the robustness of each image. Given a variance $\nu$, each $x \in T$ is replaced by a point $x'$ sampled from a Gaussian in $\mathbb{R}^{64 \times 64}$ with mean $x$ and variance $\sigma$. Denote the resulting space as $T^{\nu}$. We then define $\mathcal{S}^{\nu}_{s}$ to be a set of subsets obtained by replacing any $x \in A \in \mathcal{S}_{s}$ with $x'$. Given in Figures \ref{torus_result_0.1} and \ref{torus_result_0.3} are the results of parameter selection for $s$ (Section $5$) applied to $T^{\nu}$ with bootstrap samples $\mathcal{S}^{\nu}_s$ for both cross-matched and reference images for $\nu = 0.1,0.3$ (results for other $\nu$ are similar). It is observable that for each $\nu$, as in the noiseless case, the cross-matched and reference images agree on approximate locations of prevalent intervals. We observe with increasing $\nu$ that the lower persistent feature appears less prevalent in cross-matched images. This is apparent in the row sums directly to the right of each image: the difference in peaks is greater for cross-matched images than corresponding reference images for higher $\nu$.

\subsection{Noise Robustness of Images}
We now study the noise robustness of Algorithm \ref{alg_1} against that of the reference prevalence image. We compare $L^1$ distances between the outputs for $T$ and $T^\nu$ for both reference and cross-matched prevalence images with $s = 150$ (the $s$ value outputted following parameter selection for both schemes in the noiseless case). We also compute the $L^1$ distances between the persistence images of $T$ and $T^\nu$.

\begin{figure}
    \centering
    \includegraphics[scale = 0.45]{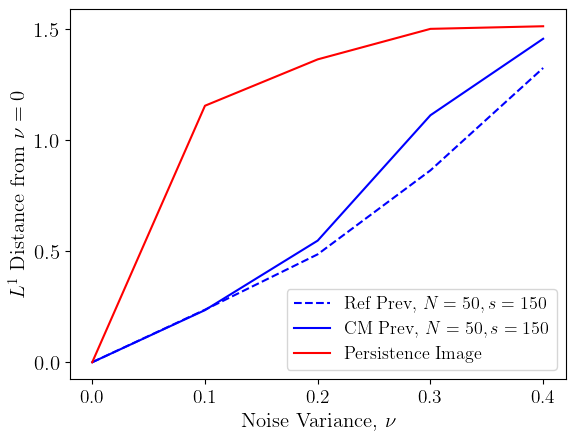}
    \caption{Noise Self-Robustness of Torus Images}
    \label{noise_robustness}
\end{figure}

In Figure \ref{noise_robustness}, we report these distances for each image across noise variances $\nu = 0.1,0.2,0.3,0.4$. We observe that the performance of both cross-matched and prevalence images are similar, both significantly outperforming the associated persistence image.

%% file: sec/7_conclusion.tex
\section{Conclusion}
In this paper, we introduced a novel method that iterates on previously proposed topological bootstrapping schemes to extract topological information from a point cloud $X$. Our method addresses shortcomings of traditional persistent homology schemes, namely (1) accuracy in the setting of $X$ nosily sampled from a low-dimensional shape embedded in high dimension and (2) poor computational complexity scaling with the size of $X$. We have demonstrated on synthetic data that this method accurately identifies prevalent topological features compared to previously proposed methods (Section \ref{sec:exp}).

\subsection{Future Work}
Our present implementation for selection of $s$ (Section \ref{sec:param}) requires the computation of potentially many cross-matched prevalence images for varying $s$. While this selection process was effective in experimentation (Section \ref{sec:exp}), it can be time intensive in practice. We have observed that standard persistent homology calculations are much faster than image persistent homology calculations for clouds of the same size. For this reason, the development of parameter selection techniques that do not require as frequent computation of image persistence may be worthwhile. Such techniques may be feasible as matching affinity is, in part, dependent on the data of sample bars (Equation \ref{eq:aff}).

%% file: cross_matching.bbl
\begin{thebibliography}{}
\setlength{\itemindent}{-\leftmargin}
\small
\bibitem{Adams} Adams, H., T.~Emerson, M.~Kirby, R.~Neville, C.~Peterson, P.~Shipman, S.~Chepushtanova, E.~Hanson, F.~Motta, L.~Ziegelmeier (2017).
	\newblock Persistence Images: A Stable Vector Representation of Persistent Homology
	\newblock Journal of Machine Learning Research 18 1-35.
	
\bibitem{Ahmed} Ahmed, M., B.~Fasy, and C.~Wenk (2014). 
	\newblock Local persistent homology based distance between maps. 
	\newblock In SIGSPATIAL. ACM, Nov. 2014

\bibitem{Bauer} Bauer, U., M.~Schmahl (2022).
	\newblock Efficient Computation of Image Persistence.
	\newblock International symposium on Computational Geometry, 2022.
	
\bibitem{PersLay} Carriere, M., F.~Chazal, Y.~Ike, T.~Lacombe, M.~Royer, Y.~Umeda (2020).
	\newblock PersLay: A Neural Network Layer for Persistence Diagrams and New Graph Topological Signatures.  					\newblock Proceedings of the Twenty Third International Conference on Artificial Intelligence and Statistics, PMLR 108:2786-2796.

\bibitem{Chen2011} Chen, C., D.~Freedman, C.H.~Lampert (2011).
\newblock Enforcing topological constraints in random field image segmentation
\newblock CVPR 2011, Colorado Springs, CO, USA, 2011, pp. 2089-2096, doi: 10.1109/CVPR.2011.5995503.

\bibitem{Deq} Dequ´eant, M., S.~Ahnert, H.~Edelsbrunner, T.~Fink, E.~Glynn, G.~Hattem,
A.~Kudlicki, Y.~Mileyko, J.~Morton, A.~Mushegian, L.~Pachter, M.~Rowicka, A.~Shiu,
B.~Sturmfels, and O.~Pourqui´e (2008)
	\newblock Comparison of pattern detection methods in microarray time series of the segmentation clock. 
	\newblock PLoS ONE, 3(8):e2856.

\bibitem{Easley} Easley, T., K.~Freese, E.~Munch, J.~Bijsterbosch (2023).
	\newblock Comparing representations of high-dimensional data with persistent homology: a case study in neuroimaging.
	\newblock arXiv 2306.13802.
	
\bibitem{ELZ} Edelsbrunner, H., D.~Letscher, and A.~Zomorodian (2000).
	\newblock Topological persistence and simplification.
	\newblock In Foundations of Computer Science. Proceedings. 41st Annual Symposium on, pages 454-463.

\bibitem{GR} Garcia-Redondo, I., A.~Monod, and A.~Song (2024).
    \newblock Fast Topological Signal Identification and Persistent Cohomological Cycle Matching.
    \newblock Journal of Applied and Computational Topology 8, 695-726.
    
\bibitem{Headd} Headd, J., Y.E.~Ban, P.~Brown, H.~Edelsbrunner, M.~Vaidya, J.~Rudolph (2007).
	\newblock Protein-protein interfaces: properties, preferences, and projections.
	\newblock Journal of Proteome Research, 6(7): 2576-2586.
	\newblock PMID: 1754262
	
\bibitem{Hensel} Hensel, F., M.~Moor, B.~Rieck (2021).
	\newblock A Survey of Topological Machine Learning Methods
	\newblock Frontiers in Artificial Intelligence vol 4.

\bibitem{Hu} Hu, X., F.~Li, D.~Samaras, and C.~Chen (2019).  
	\newblock Topology-preserving deep image segmentation.
	\newblock In Advances in Neural Information Processing Systems, 5658–5669.

\bibitem{ellipse} Jin, Y. (2022).
	\newblock ellipse (version 0.6.0)
	\newblock https://pypi.org/project/ellipse/

\bibitem{Li} Li, C., M.~Ovsjanikov, and F.~Chazal (2014). 
	\newblock Persistence-based structural recognition. 
	\newblock In Computer Vision and Pattern Recognition (CVPR), 2014 IEEE Conference on, pages 2003–2010, June 2014.

\bibitem{Myers} Myers, A., H.~Kvinge, T.~Emerson (2023).
\newblock TopFusion: Using Topological Feature Space for Fusion and Imputation in Multi-Modal Data.
\newblock Proceedings of the IEEE/CVF Conference on Computer Vision and Pattern Recognition (CVPR) Workshops, 2023, pp. 600-609.

\bibitem{Otter} Otter, N., M.A.~Porter, U.~Tillmann et al.(2017)
	\newblock A roadmap for the computation of persistent homology.
	\newblock EPJ Data Sci. 6, 17
	
\bibitem{Reani} Reani, Y., O.~Bobrowski (2023).
    \newblock Cycle Registration in Persistent Homology with Applications in Topological Bootstrap.
    \newblock IEEE Transactions on Pattern Analysis and Machine Intelligence, vol. 45, no. 5, pp. 5579-5593, 1 May 2023.

\bibitem{Seversky2016} Seversky, L., S.~Davis, and M.~Berger (2016).
\newblock On Time-series Topological Data Analysis: New Data and Opportunities.
\newblock Proceedings of the IEEE Conference on Computer Vision and Pattern Recognition (CVPR) Workshops, 2016, pp. 59-67.

\bibitem{VerHoef} Ver Hoef, L., H. Adams, E. J. King, and I. Ebert-Uphoff (2023).
\newblock A Primer on Topological Data Analysis to Support Image Analysis Tasks in Environmental Science.
\newblock Artif. Intell. Earth Syst., 2, e220039, https://doi.org/10.1175/AIES-D-22-0039.1.

\bibitem{Verri} Verri, A., C.~Uras, P.~Frosini, M.~Ferri (1993).
	\newblock On the use of size functions for shape analysis.
	\newblock Biological Cybernetics, 70, 99-107.
\end{thebibliography}
